\documentclass[letterpaper, 11pt]{article}
\usepackage{amsmath, amsthm, amssymb}
\usepackage{thm-restate}
\usepackage{graphicx}
\usepackage{xcolor}
\usepackage{algorithm,algorithmic}
\usepackage[colorlinks,linkcolor=blue,citecolor=blue,urlcolor=magenta,linktocpage,plainpages=false]{hyperref}
\usepackage{caption}
\usepackage{subcaption}
\usepackage{setspace}
\usepackage[left=1in, right=1in, top=1in, bottom=1in]{geometry}
\usepackage{authblk}
\usepackage{changepage}
\usepackage{enumitem}
\usepackage{booktabs}
\usepackage{multirow}
\newcommand{\verbrm}[1]{\rm\texttt{#1}}

\title{
Exploiting Instance and Variable Similarity to Improve Learning-Enhanced Branching
}
\author{
Xiaoyi Gu\thanks{xiaoyigu@gatech.edu, H. Milton  Stewart School of Industrial \& Systems Engineering, Georgia Institute of Technology, Atlanta, GA 30332.},
Santanu S. Dey\thanks{santanu.dey@isye.gatech.edu, H. Milton  Stewart School of Industrial \& Systems Engineering, Georgia Institute of Technology, Atlanta, GA 30332.}, 
\'Alinson S. Xavier\thanks{axavier@anl.gov, Energy Systems Division, Argonne National Laboratory}, Feng Qiu\thanks{fqiu@anl.gov, Energy Systems Division, Argonne National Laboratory}
} 

\date{Aug 21, 2022}

\begin{document}

\maketitle
\begin{abstract} 
In many operational applications, it is necessary to routinely find, within a very limited time window, provably good solutions to challenging mixed-integer linear programming (MILP) problems. An example is the Security-Constrained Unit Commitment (SCUC) problem, solved daily to clear the day-ahead electricity markets. Previous research demonstrated that machine learning (ML) methods can produce high-quality heuristic solutions to combinatorial problems, but proving the optimality of these solutions, even with recently-proposed learning-enhanced branching methods, can still be time-consuming. In this paper, we propose a simple modification to improve the performance of learning-enhanced branching methods based on the key observation that, in such operational applications, instances are significantly similar to each other. Specifically, instances typically share the same size and problem structure, with slight differences only on matrix coefficients, right-hand sides and objective function. In addition, certain groups of variables within a given instance are also typically similar to each other. Therefore, unlike previous works in the literature which predicted all branching scores with a single ML model, we propose training separate ML models per variable or per groups of variables, based on their similarity. We evaluate this enhancement on realistic large-scale SCUC instances and we obtain significantly better gap closures than previous works with the same amount of training data.
\end{abstract}
%



\section{Introduction}

Many practical optimization problems related to business or industrial \emph{operations} require the solution of challenging Mixed-Integer Programming Problems (MILP) on a regular basis. Unlike MILPs related to \emph{long-term planning}, these operational problems typically need to be solved within a limited time window, since problem data only becomes available a few hours or even minutes before the solution is to be implemented in the real world, and therefore impose great challenges to even state-of-the-art MILP solvers.
An example is the \emph{Security-Constrained Unit Commitment} (SCUC) problem, solved daily by Independent System Operators to clear the day-ahead electricity markets. Because operators only have 3 to 4 hours to clear the markets after receiving bids and offers from market participants, SCUC needs to be solved within 15 to 30 minutes~\cite{chen2016improving}. Another example is \emph{Optimal Transmission Switching} (OTS)~\cite{blumsack2007quantitative,fisher2008optimal,kocuk2016cycle,johnson2020k}, which aims to optimize the topology of the power grid. To be used in real-time operations, OTS needs to be solved every 15 minutes, within approximately 5 minutes to allow time for feasibility checks and other post-processing.

Previous research has demonstrated that machine learning (ML) methods can quickly produce high-quality heuristic solutions to a variety of challenging combinatorial problems, including some operational problems. For example, \cite{johnson2020k} uses k-nearest neighbors to predict solutions for the real-time DC optimal transmission switching; \cite{oroojlooyjadid2020applying} applies deep learning to construct solutions to the newsvendor problem; \cite{xavier2021learning} uses machine learning to construct partial solutions to the security-constraint unit commitment problem; and \cite{larsen2022predicting} uses deep learning to predict solutions to stochastic load planning problems. Proving the optimality of these heuristic solutions, however, may still require the exploration of large branch-and-bound trees.
Motivated by this, previous research has also shown that ML methods can help MILP solvers to explore the branch-and-bound tree more efficiently through \emph{learning-enhanced branching rules}. For example, \cite{alvarez2014supervised} approximates strong branching decisions using supervised learning and extremely-randomized trees; \cite{khalil2016learning} models the variable branching task as an online \emph{learning-to-rank} problem; \cite{gasse2019exact} mimics strong branching decisions using graph convolutional networks and a variable-constraint bipartite graph; \cite{nair2020solving} uses deep learning to mimic a variant of full strong branching that scales to large instances using GPUs; \cite{prouvost2020ecole} and \cite{zhang2022deep} use reinforcement learning to guide branching decisions.

In this paper, we propose a simple enhancement to improve the performance of learning-enhanced branching methods based on the key observation that, in operational problems, instances that are solved on a regular basis share significant similarities with each other.
Specifically, we routinely need to solve problems of the form
\begin{align}\label{eq:optsetup}
\begin{array}{rl}
	\textup{minimize} \quad & c_x^T x + c_y^Ty \\
	\text{subject to} \quad & Ax + Gy \geq b \\
	& x \in \{0,1\}^n \\
    & y \in \mathbb{R}^m 
\end{array}
\end{align}
where the matrices $A$ and $G$, and the vectors $c$ and $b$ only differ slightly.
The identity of each decision variable, therefore, remains the same across instances. In addition, certain groups of variables within a given instance can also be naturally grouped, based on their physical representaiton, and these groups remain the same across multiple instances.
For example, SCUC includes groups of decision variables correponding to the same generator, and groups of variables corresponding to the same operation, such as switching on a generator.
We propose, therefore, training separate ML models per variable or per group of variables, based on their meaning and similarities.
This is in contrast with previous works, which used a single ML model to predict branching decisions for all variables.
While sophisticated machine learning models could potentially learn these groupings automatically, given enough training data, we stress that generating training data for learning-enhanced branching methods can be extremely time-consuming, especially for large-scale MILP instances, as highlighted in previous works \cite{nair2020solving}. Furthermore, the groupings we propose are straightforward to identify and require negligible manual labeling effort.

To evaluate the proposed enhancement, we modify the supervised learning method proposed by \cite{alvarez2014supervised} and perform comprehensive benchmarks on realistic large-scale SCUC instances. We carefully evaluate the impact of various different groupings, including ``per-variable'', ``per-generator'', ``per-time'' and ``per-type'' on ML accuracy scores, as well as MIP gap closure. Our experiments show that variable grouping can significantly improve the performance of the original method, bringing it much closer to strong branching, while using exactly the same training data set.

The rest of the paper is organized as follows. In Section~\ref{sec:background}, we review SCUC and existing learning-enhanced branching rules. In Section~\ref{sec:proposed}, we explain our proposed branching schemes in further detail. In Section~\ref{sec:exp}, we describe our experimental setup, including implementation details, training data generation and hyperparameters. In Section~\ref{sec:evaluation} we present our experimental results, and finally in Section~\ref{sec:final} we discuss some conclusions and future directions of research.

\section{Background \label{sec:background}}

In this section, we review the Security-Constrained Unit Commitment problem, the branch-and-bound algorithm and previously-proposed learning-enhanced branching rules.

\subsection{Security-Constrained Unit Commitment \label{sec:introscuc}}
\emph{Security-Constrained Unit Commitment} (SCUC) is a fundamental optimization problem solved daily by power system operators that seeks to minimize electricity production costs by finding the most effective generator commitment schedule and power output levels and is widely used in power system applications such as electricity market clearing and reliability assessment. 
SCUC is a generalization of the \emph{unit commitment} problem (\cite{cohen1987optimization}), and is typically modeled in industrial practice as a large-scale Mixed Integer Linear Programming (MILP) problem. Numerous MILP formulations for the problem have been proposed (see, for example,~\cite{garver1962power,ostrowski2011tight,morales2013tight}).
Besides being difficult in the theoretical sense (see~\cite{bendotti2019complexity} for NP-hardness of SCUC), SCUC is made even more challenging by the fact that a new near-optimal solution must be obtained every day within a very limited time window.

\subsection{Branch-and-Bound Algorithm}\label{sec:introbranch}
The branch-and-bound algorithm was proposed in the landmark paper \cite{Land60anautomatic} to provide a finite algorithm to solve MILPs, and today it is at the heart of all state-of-the-art MILP solvers. To completely define the algorithm, it is necessary to specify a \emph{variable selection rule}, which decides how to partition the problem into multiple subproblems, and a \emph{node selection rule}, which defines what subproblem to process next (\cite{wolsey1999integer,conforti2014integer}). On the node selection side, it is well-known that using the \textit{best-bound rule} produces the smallest trees for realistic benchmark instances (\cite{achterberg2007constraint}). Recently, \cite{dey2020branch} also showed that, for certain classes of random integer programs, using this rule guarantees a polynomial-size branch-and-bound tree with high probability.
The more challenging rule to optimize is the variable selection rule. Negative results due to~\cite{jeroslow1974trivial,chvatal1980hard,cook1990complexity, dadush2020complexity,jiang2021complexity,dey2022lower,dey2021lower} show that, no matter what variable selection rule is used, the branch-and-bound tree is, of exponential size in the worst case. In practice, the variable selection rule has a huge influence on the size of the branch-and-bound tree. Today, it is well-established that \textit{strong branching} (\cite{applegate1995finding}) typically produces the smallest trees for a variety of optimization problems. Although this rule produces high-quality decisions (\cite{dey2021theoretical}), it is extremely computationally demanding, and most times prohibitive to implement (\cite{achterberg2005branching}). 
See~\cite{linderoth1999computational,lodireview} for a nice overview of computational issues with the branch-and-bound algorithm. 
Since strong branching produces high-quality decisions, but at a high computational cost, multiple works in the literature have attempted to mimic it using machine learning methods~\cite{khalil2016learning,alvarez2017machine,gasse2019exact,yang2020learning,gupta2020hybrid,nair2020solving}. See~\cite{lodi2017learning}  for a review of this direction of research.

\subsection{Non-ML Branching Methods}
Most traditional non-ML branching methods can be generally modelled as a score function $S(i,I)$, where $i$ represents the index of the candidate variable and $I$ represents general information about current and previous nodes. Given this function, the candidate variable chosen for branching is the one with the best (highest) score.

\textit{Most-infeasible branching} (MIB) is one of the most simple branching strategies, which uses a score function of $S_{MIB}(i, I) = \min\{x_i,1-x_i\}$, the infeasibility or ``fraction'' of the candidate variable at current node. 
While \verb|MIB| is extremely cheap to evaluate, it is unfortunately known to perform poorly compared to other branching methods.
\textit{Reliability branching} (\verbrm{RB:$\lambda$:$\eta$}) is the state-of-the-art branching scheme with respect to running time (\cite{achterberg2005branching}).
The method relies on \textit{probing}, during which it actually solves the linear relaxations of both upward and downward branches for the fractional variable. More specifically, at each node, the candidate fractional variables are sorted by pseudocost, then at most $\lambda$ variables are probed. If a given variable has already been probed $\eta$ times during the entire execution of the algorithm, its pseudocosts are deemed \emph{reliable}, and no more probes are performed for the variable. The score function is usually $S(i, I) = \max\{\tilde{\Delta}_i^-,\epsilon\}\cdot\max\{\tilde{\Delta}_i^+,\epsilon\}$, where $\tilde{\Delta}_i^\pm$ are the objective increases if probed, or the pseudocost estimations otherwise.

The scheme covers \verbrm{RB:inf:inf}, so called \textit{full strong branching}, where \textit{every} fractional variable is probed, and which was shown in numerous results (like \cite{achterberg2005branching}) to generate the smallest branch-and-bound trees compared to other rules. However, since every fractional variable to be probed translates to two linear relaxations to be solved, this scheme is computationally prohibitive for large scale problems.  
On the other hand, reduced versions of \verbrm{RB:$\lambda$:$\eta$} with smaller values of $\lambda$ and $\eta$ can achieve good branching while being considerably faster.
Considering the size of our problems, we use settings of \verbrm{RB:100:inf} 
as the \emph{practical strong branching oracle}.

\subsection{Previous ML Branching Methods}
In this work, we evaluate our proposed ML scheme against the method described in \cite{alvarez2017machine}, denoted as \verb|ML:ET|. 
In general, \verb|ML:ET| can be formulated as
$$
S_{ML:ET}(i, I) = f_{ML:ET}(\phi_i)\approx S_{oracle}(i,I),
$$
where $\phi_i$ is the feature vector from available information $I$, and $f_{ML:ET}$ is the learned machine-learning model, which is used to mimic a good branching oracle like strong branching.
A wide variety of handcrafted features, including static problem features, dynamic problem features and dynamic optimization features, are used.
Training problem instances are solved using strong branching, without any heuristics or cut generation, to collect strong branching scores, and then Extremely Randomized Trees (ExtraTrees, \cite{geurts2006extremely}) are used to approximate them.

\section{Proposed Per-Variable and Per-Group ML Methods \label{sec:proposed}}
The \emph{per-variable ML scheme} we propose, denoted by
\verbrm{ML:PV},
looks similar to \verbrm{ML:ET}, with the key difference that one model is built for each variable.
In general, our scheme can be formulated as
$$
S_{ML:PV}(i,I) = f^i_{ML:PV}(\phi_i)\approx S_{oracle}(i, I).
$$
While the idea applies to any regressor, we use ExtraTrees to keep it consistent with \verbrm{ML:ET}. However, since the size of the training dataset for each per-variable model is significantly smaller, the hyperparameters need to be modified to avoid over-fitting. Specifically, smaller trees need to be generated.

We also propose a \emph{per-group ML scheme}, denoted by
\verbrm{ML:PG},
which is similar to the per-variable scheme, which trains one model per variable group. It is formulated as
$$
S_{ML:PG}(i,I) = f^{g(i)}_{ML:PG}(\phi_i)\approx S_{oracle}(i, I),
$$
where $g(i)$ is the corresponding group to which the $i$-th column belongs. Similarly, we use ExtraTrees as the regressor and we tune the hyperparameters individually, for each group of variables, to avoid over-fitting.

Note that both \verbrm{ML:ET} and 
\verbrm{ML:PV} 
could be viewed as extreme cases of per-group settings, since \verbrm{ML:ET} is basically grouping all variables into one group, while 
\verbrm{ML:PV} 
is grouping every variable into its own group.
We anticipate more accurate predictions with finer grouping, at the cost of reduced model generality.

\section{Experimental Setup \label{sec:exp}}
In this section, we describe our branch-and-bound implementation, computational environment, benchmark instances and training data generation.

\subsection{Implementation}

To ensure consistent training data and reproducible experimental results, in this work we implemented, in the Julia programming language, a textbook version of the branch-and-bound algorithm, as well multiple variable branching rules and multiple node selection rules. While cutting planes can significantly enhance the performance of textbook branch-and-bound methods, in this work we did not consider their usage, as they make the analysis of branching rules significantly more difficult. The algorithm was also given the primal optimal value to eliminate the influence of primal heuristics. The implementation relies on an external LP solver, accessed through \texttt{MathOptInterface}, to process node relaxations and evaluate strong branching decisions. In our experiments, we used \emph{best-bound with plunging} for node selection and Gurobi 9.0 (\cite{gurobi}) as the LP solver. The branch-and-bound implementation has been made available as part of the open-source package MIPLearn (\cite{xavier2020miplearn}). MIPLearn was also used to compute variable features. Machine learning models were implemented in Python 3.9 with \verbrm{scikit-learn}, and \verbrm{PyCall.jl} was used to call the ML models from Julia.

\subsection{Environment}

Both training and testing were processed in parallel on the high-performance computational cluster of ISyE, Georgia Tech, which contains roughly 2,340 cores of x86-64 processing with over 28.9TB of memory spread across the systems and each task was processed on a dedicated core with 8GB of memory.
A time limit of one week (604,800 seconds) was imposed when solving all MILPs, though never activated, since all runs finished within the allotted time. In all runs, the branch-and-bound algorithm was also configured to terminate when a 0.01\% relative gap was reached.

\subsection{Instances}

Five realistic large-scale security-constrained unit commitment instances from \texttt{UnitCommitment.jl} \cite{xavier2020unitcommitment}, corresponding to five European transmission networks from MATPOWER (\cite{zimmerman2010matpower,fliscounakis2013contingency,josz2016ac}), were used to evaluate the performance of the ML branching rules.
For each network, 50 instance variations were generated by applying the randomization algorithm described in \cite{xavier2021learning} -- thus giving us 250 instance variations in total. For each network, the first 40 instances variations were used for training, while the remaining 10 were used for testing.
The sizes of the resulting instances are displayed in Table~\ref{tab:size}.
We highlight that the instances used in our experiments are quite large, with up to 400,000 decision variables and 350,000 constraints. To the best of our knowledge, most papers in the learning-to-branch literature have not dealt with such large instances, with the notable exception of~\cite{nair2020solving}.

\begin{table}[htbp]
  \centering
  \caption{Size of Networks}
    \begin{tabular}{c|cccc|ccc}
    \toprule
    Network & Hours & Generators & Buses & Lines & Variables & Rows  & Binaries \\ \midrule
    case1888rte & 24    & 296   & 1,888  & 2,531  & 235,591 & 196,783 & 41,232 \\
    case1951rte & 24    & 390   & 1,951  & 2,596  & 266,088 & 244,220 & 54,144 \\
    case2848rte & 24    & 544   & 2,848  & 3,776  & 377,760 & 340,228 & 71,904 \\
    case3012wp & 24    & 496   & 3,012  & 3,572  & 357,146 & 305,076 & 59,712 \\
    case3375wp & 24    & 590   & 3,374  & 4,161  & 413,161 & 357,065 & 71,856 \\
    \bottomrule
    \end{tabular}%
    \label{tab:size}%
\end{table}%

The package \texttt{UnitCommitment.jl} was also used to construct the MILP, using a state-of-the-art formulation of the problem. In the following, we list the description of all binary variables in the formulated MILP, as it is of particular interest for both per-variable and per-group approach. We refer to the package documentation for more details: 
\begin{itemize}
    \item \verbrm{is\_on[g,t]}: True if generator $g$ is on at time $t$
    \item \verbrm{switch\_on[g,t]}: True if generator $g$ switches on at time $t$.
    \item 
    \verbrm{switch\_off[g,t]}: True if generator $g$ switches off at time $t$.
    \item 
    \verbrm{startup[g,t,s]}: True if generator $g$ switches on at time $t$ incurring start-up costs from start-up category $s$.
\end{itemize}

\subsection{Training dataset}

To generate the training dataset,  we solved each training instance using \verbrm{RB:100:inf} and, for each strong branching call made by the algorithm, we collected features describing the evaluated variable, as well as the logarithm of the computed strong branching score (due to the fact that it spans several orders of magnitude).
We enforced a node limit of 1000, a time limit of $1.2\times 10^6$ seconds (roughly 2 weeks, due to server limitation) and a relative gap limit of 0.01\%. The algorithm stopped when any of these limits were reached.

\subsection{ML models \label{subsec:mlmodels}}

For each network, one model was trained for \verb|ML:ET|, while a series of models were trained for the per-variable scheme or the per-group schemes.
To the best of our ability and understanding, \verb|ML:ET| is an accurate implementation of the method described in \cite{alvarez2017machine}. In particular, we used the same set of features and the same ML regressors. We also performed hyperparameters tuning, using Scikit-Learn's \verb|GridSearchCV|, to improve the classifier's performance and to manage its memory requirements. Our final set of hyperparameters for \verb|ML:ET| were \verb|min_samples_split=10|, \verb|max_depth=25| and \verb|n_estimators=25|.

For the per-group scheme, we chose three different grouping methods: 
\begin{itemize}
    \item \emph{Per-generator} (\verbrm{ML:PGE}): Two variables are grouped together if they have the same base name (e.g. \verbrm{is\_on}, \verbrm{switch\_on}, \verbrm{switch\_off}), the same generator index $g$ and the same startup category $s$. Each group, therefore, consists of 24 variables, corresponding to different values of time $t$.

    \item \emph{Per-time} (\verbrm{ML:PTI}): Variables are grouped together if they have the same base name, the same time $t$ and the same startup category $s$. We allow, therefore, grouping of variables corresponding to different generators.

    \item \emph{Per-name} (\verbrm{ML:PNA}): Variables are grouped together if the have the same base name and the same startup category $s$. Since the benchmark instances have three different startup categories, we have exactly 6 groups in total. This grouping strategy was the closest to \verbrm{ML:ET} in the number of models.
\end{itemize}

As described in Section~\ref{sec:proposed}, we also used ExtraTrees as classifier for per-group strategies, but a different set of hyperparameters, found through Scikit-Learn's \verb|GridSearchCV|, to avoid over-fitting. For all regressors, we used \verb|n_estimators=25|. For per-variable \verb|ML:PV| and the per-generator \verbrm{ML:PGE}, we selected \verb|min_samples_split=5| and \verb|max_depth=12|. For per-time \verb|ML:PTI|, we used \verb|min_samples_split=8| and \verb|max_depth=12|. For per-name \verb|ML:PNA|, we used \verb|min_samples_split=8| and \verb|max_depth=16|.

In Table~\ref{tab:nummodel}, we display the number of groups for each grouping strategy and for each network. Note that these numbers indicate the \emph{maximum} number of ML models trained; if the strong branching routine is never called during training for all variables in certain group, or if it is only called less than 10 times, then not enough training data is available for that group, and therefore no ML model is trained. We recall that the strong branching routine is not called for a given variable if either the variable is already integral at the node, or if the maximum number of strong branching evaluations has been reached. The actual number of generated models is displayed in Table~\ref{tab:actualnum}. We noticed that the difference between the number of generated models and the theoretical limits gets more prominent as the groups become finer, which is expected. During inference time, if a per-variable or per-group model is not available, we fall back to the general-purpose \verbrm{ML:ET} model.

\begin{table}[htbp]
  \centering
  \caption{Maximum number of ML models for different grouping strategies.}
    \begin{tabular}{cccccc}
      \toprule
    Network & ML:ET & ML:PNA & ML:PTI & ML:PGE & ML:PV \\\midrule
    case1888rte     & 1     & 6     & 144   & 1,718  & 41,232 \\
    case1951rte     & 1     & 6     & 144   & 2,256  & 54,144 \\
    case2848rte     & 1     & 6     & 144   & 2,996  & 71,904 \\
    case3012wp     & 1     & 6     & 144   & 2,488  & 59,712 \\
    case3375wp     & 1     & 6     & 144   & 2,994  & 71,856 \\\bottomrule
    \end{tabular}%
  \label{tab:nummodel}%
\end{table}%

\begin{table}[htbp]
  \centering
  \caption{Actual number of ML models for different grouping strategies.}
    \begin{tabular}{cccccc}
      \toprule
    Network & ML:ET & ML:PNA & ML:PTI & ML:PGE & ML:PV \\\midrule
    case1888rte     & 1     & 6     & 130   & 488   & 3,156 \\
    case1951rte     & 1     & 6     & 123   & 406   & 2,957 \\
    case2848rte     & 1     & 6     & 129   & 502   & 3,076 \\
    case3012wp     & 1     & 6     & 127   & 479   & 3,053 \\
    case3375wp     & 1     & 6     & 126   & 703   & 3,710 \\\bottomrule
    \end{tabular}%
  \label{tab:actualnum}%
\end{table}%

\section{Evaluation}\label{sec:evaluation}
In this section, we evaluate the performance of different grouping strategies on SCUC problems. In Subsection~\ref{subsec:accuracy}, we start by evaluating the accuracy of different ML models in predicting strong branching scores. In Subsections~\ref{subsec:smalltree} and \ref{subsec:largetree}, we evaluate their impact on MIP gap closure after 1,000 and 10,000 branch-and-bound nodes, respectively. Finally, in Subsection~\ref{subsec:presolve}, we evaluate the impact of presolve on the effectiveness of different methods.

\subsection{Model accuracy \label{subsec:accuracy}}

To measure the accuracy of different grouping strategies, we used $k$-fold cross-validation, a standard approach to test the quality of ML classifiers and regressors. The training data was split into $k$ parts (folds), then $k$ estimators were trained on data from $k-1$ folds and their performance was evaluated on the remaining fold. See~\cite{hastie2009elements} for details.

\begin{figure}
     \centering
     \begin{subfigure}[b]{0.325\textwidth}
    \centering
    \includegraphics[width=\textwidth]{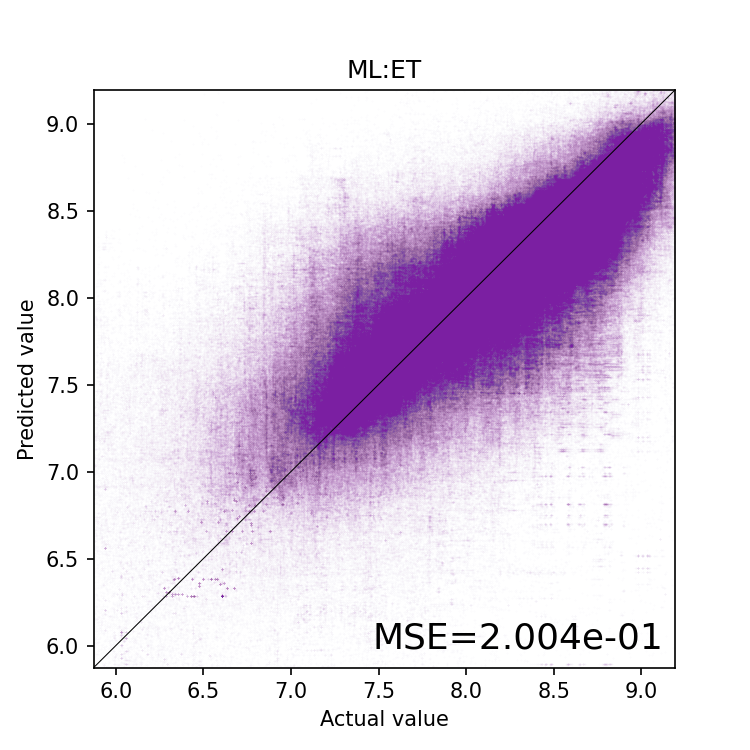}
     \end{subfigure}
     \begin{subfigure}[b]{0.325\textwidth}
    \centering
    \includegraphics[width=\textwidth]{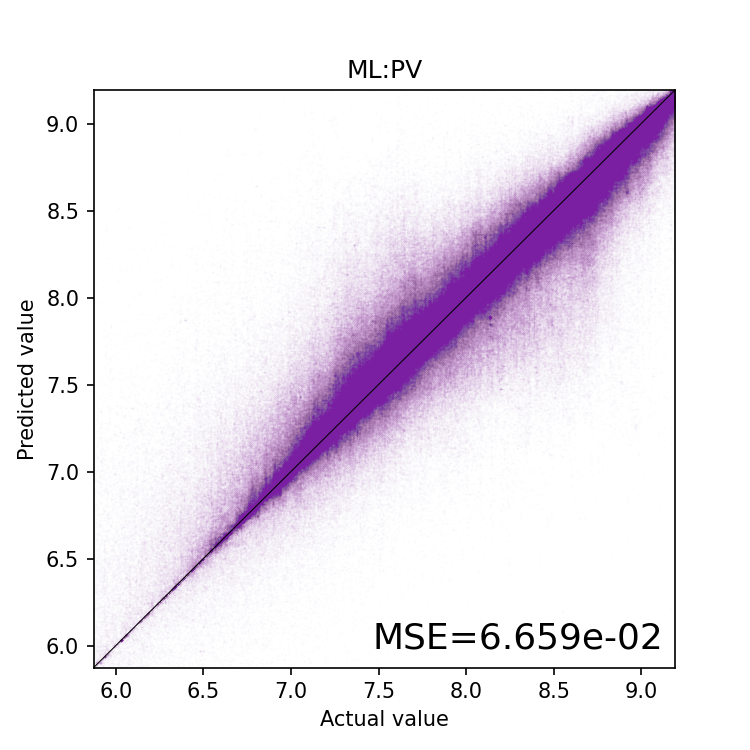}
     \end{subfigure}
     
    \centering
     \begin{subfigure}[b]{0.325\textwidth}
    \centering
    \includegraphics[width=\textwidth]{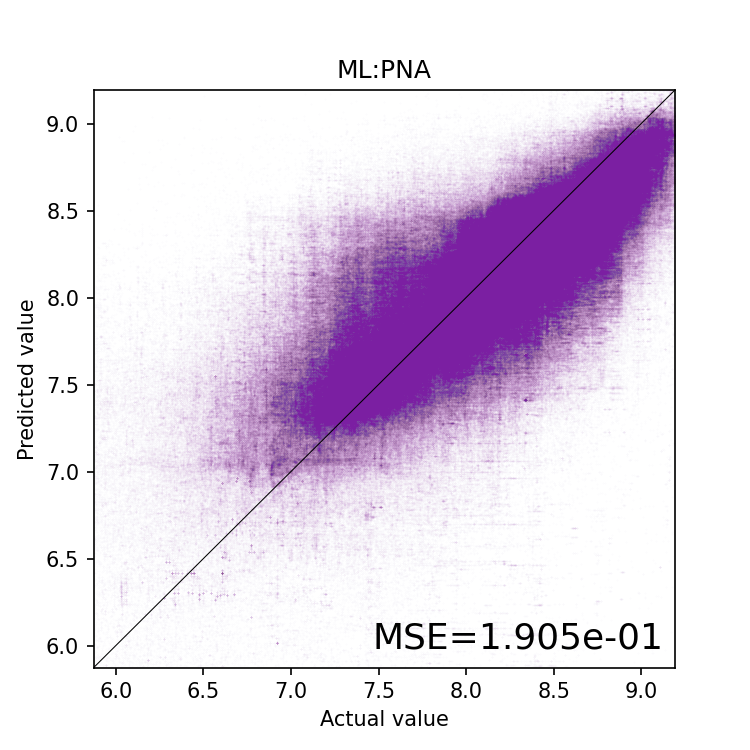}
     \end{subfigure}
     \begin{subfigure}[b]{0.325\textwidth}
    \centering
    \includegraphics[width=\textwidth]{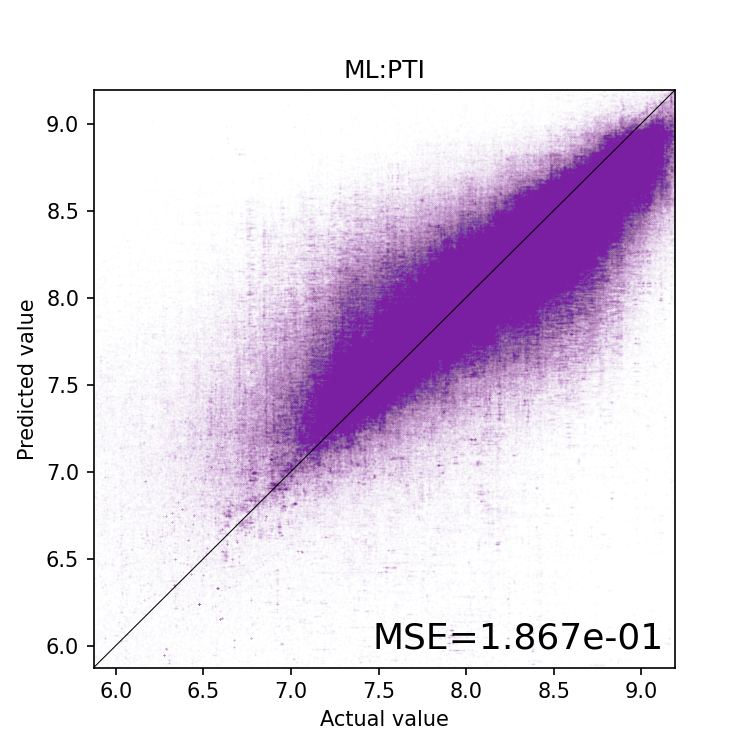}
     \end{subfigure}
     \begin{subfigure}[b]{0.325\textwidth}
    \centering
    \includegraphics[width=\textwidth]{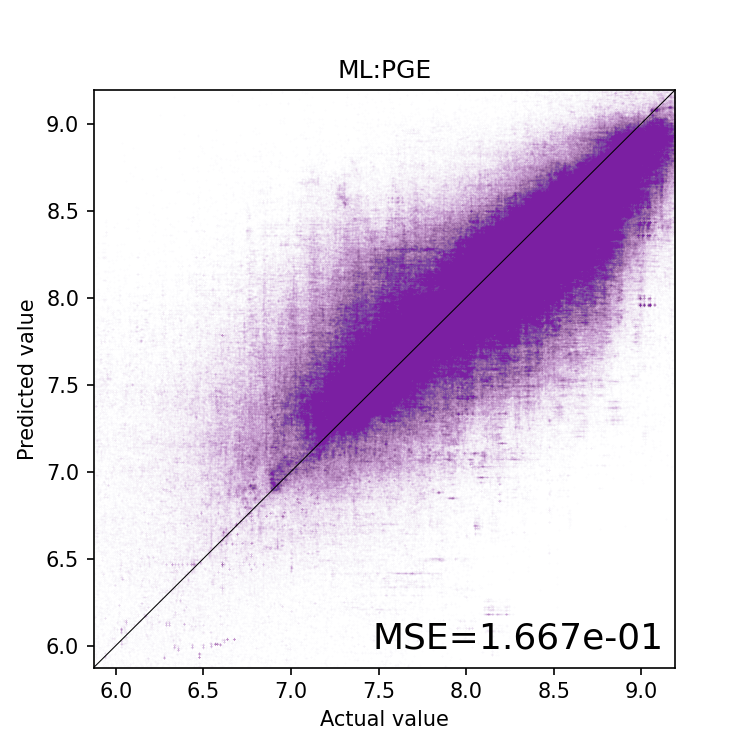}
     \end{subfigure}     
     \caption{Comparison of actual versus predicted strong branching scores on network case1888rte.}
     \label{fig:cvscore}
\end{figure}

Figure~\ref{fig:cvscore} shows the 5-fold cross-validated \emph{mean squared error} (MSE) for all grouping strategies on one particular network (case1888rte). We omit other networks, since they presented similar results.  In the charts, ``actual value'' corresponds to the strong branching scores, as generated by the non-ML branching oracle \verbrm{RB:100:inf} during training data collection, while ``predicted value'' is the branching score predicted by each ML model.

While all grouping strategies were able to approximate strong branching score to some extent, MSE scores got better as grouping became finer and each group became smaller. This was especially true for the per-variable model \verbrm{ML:PV}, the finest grouping, which achieved much better MSE results than all other methods, including \verbrm{ML:ET}.

\subsection{Gap closure on small trees (1,000 node limit) \label{subsec:smalltree}}

In the previous subsection, we evaluated the performance of the ML models in isolation, based on MSE. In this subsection, we integrate the models into the branch-and-bound algorithm and evaluate their effectiveness at solving an actual MILP.
Our main focus is the quality of branching,
measured by the final relative MIP gap after a certain number of branch-and-bound nodes.
We do not focus on total running time since this metric depends heavily on the quality of the software implementation, instead of simply on the fundamental properties of the branching rule.
It is also self-evident, in our opinion, that all evaluated ML-based branching rules, if carefully implemented, would require significantly less computational effort per node than strong branching, since they involve only the evaluation of a few small decision trees per variable, instead of the solution of a large-scale linear programming problem.

\begin{table}[t]
  \centering
  \caption{Relative MIP gap after at most 1,000 branch-and-bound nodes.}
    \begin{tabular}{crrrrrrr}
    \toprule
    \multirow{2}{*}{Network} & \multicolumn{7}{c}{Relative MIP gap (\%)} \\
    \cline{2-8}
    & \tt MIB & \tt RB:100:inf & \tt ML:ET & \tt ML:PNA & \tt ML:PTI & \tt ML:PGE & \tt ML:PV \\\midrule
    case1888rte & 1.78  & 0.76  & 1.31  & 1.35  & 0.90  & 1.22  & \bf 0.88 \\
    case1951rte & 0.41  & 0.19  & 0.23  & 0.24  & 0.20  & 0.22  & \bf 0.20 \\
    case2848rte & 0.83  & 0.37  & 0.54  & 0.59  & 0.45  & 0.58  & \bf 0.41 \\
    case3012wp & 0.29  & 0.10   & 0.08  & \bf 0.01  & 0.02  & \bf 0.01  & 0.04 \\
    case3375wp & 0.48  & 0.14  & 0.52  & 0.46  & 0.50  & \bf 0.41  & 0.48 \\\midrule
    Average & 0.76  & 0.31  & 0.54  & 0.53  & 0.41  & 0.49  & \bf 0.40 \\
    \bottomrule
    \end{tabular}%
  \label{tab:gap1000}%
\end{table}%

In our first set of experiments, we compared the relative MIP gap attained by seven different ML and non-ML branching methods (\verbrm{MIB}, \verbrm{RB:100:inf}, \verbrm{ML:ET}, \verbrm{ML:PNA}, \verbrm{ML:PTI}, \verbrm{ML:PGE} and \verbrm{ML:PV}) after the exploration of at most 1,000 branch-and-bound nodes.
This node limit was set relatively low because we wanted to include \verbrm{RB:100:inf}, a computationally expensive branching rule.
Note that even though \verbrm{RB:100:inf} is faster than full strong branching (\verbrm{RB:inf:inf}), for the size of our instances, it is not fast enough to be a practical branching method. 
Table~\ref{tab:gap1000} presents a summary of the results. We remind the reader that 10 test instances were solved for each network; therefore, each cell in the table represents the average of 10 values. For each network, we also highlighted the ML branching rule with best performance.

As clearly illustrated, the per-variable approach \verbrm{ML:PV} outperformed \verbrm{ML:ET} on almost all networks. Moreover, \verbrm{ML:PV} presented a relative MIP gap of 0.40\%, which is only about 22\% worse than \verbrm{RB:100:inf}, and therefore provides a relatively close approximation of strong branching.

Per-generator \verbrm{ML:PGE} and per-time \verbrm{ML:PTI} approaches also significantly outperformed \verbrm{ML:ET}, although to a smaller degree in most instances. The per-name \verbrm{ML:PNA} approach overall offered little improvement, likely due to its similarity to the original \verbrm{ML:ET}.

As explained in Subsection~\ref{subsec:mlmodels}, when we need to estimate the strong branching score of a variable that does not have a per-group ML model, we fall back to the most general \verbrm{ML:ET} model. It is thus worth investigating how often such phenomenon occurred. Table~\ref{tab:fallbackrate} shows what percentage of all evaluations relied on the fallback model, instead of the specific per-group model, due to missing data. 

It is clear that, with finer grouping strategies, fallback happened more often, indicating that fewer variables were covered during training. This highlights the need to balance between more groups and more data per group. While finer groups leads to improved accuracy, this comes at the cost of reduced coverage. Nevertheless, it is also clear that the fallback rate is quite low, and that fallback does not hinder the potential of the per-group or per-variable approach. With a more substantial dataset, we also expect that we should be able to reduce the occurrence of fallback or even eliminate it altogether.

\begin{table}[t]
  \centering
  \caption{Fallback Rate, \verbrm{node\_limit=1000}}
    \begin{tabular}{crrrrr}
    \toprule
    \multirow{2}{*}{Network} & \multicolumn{5}{c}{Fallback rate (\%)} \\
    \cline{2-6}
    & \tt ML:ET & \tt ML:PNA & \tt ML:PTI & \tt ML:PGE & \tt ML:PV \\
    \midrule
    case1888rte & 0.0     & 0.0     & 0.0 & 0.1 & 9.8 \\
    case1951rte & 0.0     & 0.0     & 0.2 & 0.1 & 4.5 \\
    case2848rte & 0.0     & 0.0     & 0.0 & 0.1 & 7.2 \\
    case3012wp  & 0.0     & 0.0     & 0.0 & 2.8 & 14.6 \\
    case3375wp  & 0.0     & 0.0     & 0.0 & 0.4 & 8.5 \\
    \midrule
    Average     & 0.0     & 0.0     & 0.0 & 0.7 & 9.0 \\
    \bottomrule
    \end{tabular}%
  \label{tab:fallbackrate}%
\end{table}%

\subsection{Gap closure on large trees (10,000 node limit) \label{subsec:largetree}}
While the experiments in Section~\ref{subsec:smalltree} clearly demonstrate the superiority of the per-group methods in smaller branch-and-bound trees, it is natural to ask whether these results still hold for larger trees.
To answer this question, in this section we repeat the experiments of Section~\ref{subsec:smalltree} with a larger budget of 10,000 nodes. Due to the significant computational cost of \verbrm{RB:100:inf}, we omit this method in the following experiments.

Table~\ref{tab:gap10000} shows a summary of our results. With larger trees, the performance improvement of \verbrm{ML:PV} over \verbrm{ML:ET} becomes more clear, and the per-variable method now outperforms \verbrm{ML:ET} in every network. Consistently with our previous experients, other per-group schemes, especially per-time (\verbrm{ML:PTI}) or per-generator (\verbrm{ML:PGE}), also outperform \verbrm{ML:ET} in most cases. For network case3375wp, in particular, \verbrm{ML:PGE} presents much better performance than all the other ML methods.

\begin{table}[htbp]
  \centering
  \caption{Relative MIP Gap, \verbrm{node\_limit=10000}}
    \begin{tabular}{crrrrrr}
    \toprule
    \multirow{2}{*}{Network} & \multicolumn{6}{c}{Relative MIP gap (\%)} \\
    \cline{2-7}
    & \tt MIB   & \tt ML:ET & \tt ML:PNA & \tt ML:PTI & \tt ML:PGE & \tt ML:PV \\
    \midrule
    case1888rte & 1.75  & 1.13  & 1.20  & 0.66  & 1.03  & \bf 0.61 \\
    case1951rte & 0.37  & 0.11  & 0.11  & 0.08  & 0.11  & \bf 0.08 \\
    case2848rte & 0.77  & 0.43  & 0.47  & 0.33  & 0.46  & \bf 0.30 \\
    case3012wp & 0.27  & 0.05  & \bf 0.01  & \bf 0.01  & \bf 0.01  & 0.02 \\
    case3375wp & 0.45  & 0.51  & 0.42  & 0.47  & \bf 0.31  & 0.44 \\\hline
    Average & 0.72  & 0.45  & 0.44  & 0.31  & 0.39  & \bf 0.29 \\
    \bottomrule
    \end{tabular}%
  \label{tab:gap10000}%
\end{table}%

\subsection{Presolved Instances \label{subsec:presolve}}

The approach we propose in this work fundamentally depends on variables and variable groups keeping their identity across multiple instances. While it is true that, in operational problems, one often has to solve instances where changes are only made to the matrix coefficients, objective function and right-hand side, in practice these MILPs would still go through an extra \emph{presolve step} just before being solved, where the MILP solver makes further changes to the problem in an attempt to make it easier, based on the actual problem data. These changes could potentially modify the problem structure. To determine whether the proposed scheme would still work under this scenario, we ran further computational experiments where the test instances are presolved by Gurobi prior to being solved by our branch-and-bound implementation. No presolve, however, is applied to the training instances. Because of this, during inference, we made the decision to provide to the ML models static variable features (e.g. objective coefficient) corresponding to the original problem, not the presolved one, together with dynamic features extracted at the current node (e.g. depth of the node).

Tables~\ref{tab:pregap1000} and \ref{tab:pregap10000} show the relative MIP gap results after 1,000 and 10,000 branch-and-bound nodes, respectively. As shown in the tables, even though the test instances are slightly different the training ones, the machine learning methods still presented good relative MIP gap closures, demonstrating the robustness of the method. Moreover, the per-variable and per-group approaches again outperformed \verbrm{ML:ET} on most networks, for both node limit settings.

\begin{table}[t]
  \centering
  \caption{Relative MIP Gap after 1,000 nodes, with presolved test instances.}
    \begin{tabular}{crrrrrrr}
    \toprule
    \multirow{2}{*}{Network} & \multicolumn{7}{c}{Relative MIP gap (\%)} \\
    \cline{2-8}
      & \tt MIB & \tt RB:100:inf & \tt ML:ET & \tt ML:PNA & \tt ML:PTI & \tt ML:PGE & \tt ML:PV \\
     \midrule
    case1888rte & 1.00  & 0.23  & 0.51  & 0.52  & 0.38  & 0.41  & \bf 0.33 \\
    case1951rte & 0.29  & 0.14  & 0.17  & 0.17  & \bf 0.15  & 0.17  & 0.16 \\
    case2848rte & 0.43  & 0.19  & 0.30  & 0.33  & 0.28  & 0.29  & \bf 0.26 \\
    case3012wp & 0.05  & 0.01  & \bf 0.01  & \bf 0.01  & \bf 0.01  & \bf 0.01  & \bf 0.01 \\
    case3375wp & 0.31  & 0.02  & 0.31  & 0.27  & 0.32  & \bf 0.20  & 0.31 \\\hline
    Average & 0.41  & 0.12  & 0.26  & 0.26  & 0.23  & \bf 0.22  & \bf 0.22 \\
    \bottomrule
    \end{tabular}%
  \label{tab:pregap1000}%
\end{table}%

\begin{table}[t]
  \centering
  \caption{Relative MIP Gap after 10,000 nodes, with presolved test instances.}
  \begin{tabular}{crrrrrr}
    \toprule
    \multirow{2}{*}{Network} & \multicolumn{6}{c}{Relative MIP gap (\%)} \\
    \cline{2-7}
      & \tt MIB & \tt ML:ET & \tt ML:PNA & \tt ML:PTI & \tt ML:PGE & \tt ML:PV \\
     \midrule
    case1888rte & 0.95  & 0.37  & 0.36  & 0.29  & 0.27  & \bf 0.25 \\
    case1951rte & 0.25  & \bf 0.04  & \bf 0.04  & \bf 0.04  & 0.05  & 0.05 \\
    case2848rte & 0.39  & 0.21  & 0.24  & 0.19  & 0.19  & \bf 0.17 \\
    case3012wp & 0.03  & \bf 0.01  & \bf 0.01  & \bf 0.01  & \bf 0.01  & \bf 0.01 \\
    case3375wp & 0.30  & 0.27  & 0.22  & 0.29  & \bf 0.13  & 0.25 \\\hline
    Average & 0.38  & 0.18  & 0.18  & 0.16  & \bf 0.13  & 0.14 \\
    \bottomrule
    \end{tabular}%
  \label{tab:pregap10000}%
\end{table}%

\section{Discussion and future research directions}\label{sec:final}

In this work, we have provided evidence that building separate ML models for individual variables, or groups of variables, can result in significantly better branching decisions in the context of SCUC.
Specifically, in our experiments, the proposed per-variable \verbrm{ML:PV} method outperformed the previously described \verbrm{ML:ET} in almost every instance, under various settings. We also found that \verbrm{ML:PV} could mimic strong branching decisions very well. The per-generator \verbrm{ML:PGE} and per-time \verbrm{ML:PTI} methods presented strong performance on selected networks, frequently outperforming \verbrm{ML:ET}. Both per-variable and per-group schemes offered stable and robust performance, even on presolved instances.

Finally, we discuss some future research directions.

\begin{enumerate}

    \item \textbf{Application to other ML branching methods.} While we focused on modifying the method proposed by \cite{alvarez2017machine}, note that the method can be applied to other ML branching methods as well. For example, it would be straighforward to adapt it to the online learning method proposed by \cite{marcos2016online}, and we could apply it to the learning-to-rank approach of \cite{khalil2016learning} by focusing on ranking variables within each group.  Another question is whether graph-based problem representations, as proposed by \cite{gasse2019exact}, would eliminate the advantages of variable grouping.
    In our understanding, the main advantage of such representation is handling instances with significantly different structure. When the problems have a fixed structure, which is the setting of our present work, the advantage of graph-based representation is much less clear, since the graph would remain fixed.
    We note that many other works in the learning-to-branch literature are orthogonal to our work. For example, the GPU-based strong branching approach in \cite{nair2020solving} could be applied in our work to accelerate training data generation, but, in our understanding, would not fundamentally change the results presented.
    \item \textbf{Training on small instances then testing on larger instances.} With appropriate variable grouping, it may be possible to collect data from smaller instances, where it becomes feasible to use high-quality oracles, such as \emph{full strong branching}, then use the trained ML models to solve instances of much larger scale. In the SCUC problem, for example, if the variable groups don't take \emph{time} into consideration, then training on instances that have a shorter time horizon may be a feasible strategy.

    \item \textbf{Integrating with other MILP solvers techniques.} One particularly interesting topic is integrating this method with \textit{primal heuristics} and \textit{cutting planes}, which are widely used in the current commercial solvers to accelerate solving MILP. In this work we demonstrated that the proposed per-variable and per-group approaches are robust against presolve, but it would be interesting to evaluate their performance on the presence of other solver features.

    \item \textbf{Evaluating on other operational problems.}
    In this work, we focused exclusively on SCUC, but the method could potentially benefit other operational problems. As mentioned in the introduction, another example would be the \emph{Optimal Transmission Switching problem} (OTS), which is not currently used in real-time operations due to slow computational performance.
    An important future direction is to explore if OTS and other important operational problems could also benefit from the proposed method.
\end{enumerate}

\section{Acknowledgments}

This material is based upon work supported by the U.S. Department of Energy Advanced Grid Modeling Program. Santanu S. Dey gratefully acknowledges the support by Airforce Office of Scientific Research.

\bibliographystyle{plain}
\bibliography{mybibliography}
\end{document}